\documentclass[12pt]{article}
\usepackage{graphicx}

\usepackage{amssymb,amsmath,amsthm}

\def\0b{{\bold 0}}

\def\0b{{\bold 0}}

\def\a{{\bf a}}

\def\RR{{\mathbb R}}
\def\ZZ{{\mathbb Z}}
\def\NN{{\mathbb Z}_{\geq 0}}
\def\QQ{{\mathbb Q}}

\def\e{{\bf e}}
\def\A{{A}}

\newtheorem{Theorem}{Theorem}%[section]
\newtheorem{Lemma}[Theorem]{Lemma}
\newtheorem{Corollary}[Theorem]{Corollary}
\newtheorem{Proposition}[Theorem]{Proposition}
\newtheorem{Remark}[Theorem]{Remark}

\newtheorem{Problem}[Theorem]{Problem}

\title{Non-very ample configurations arising from
contingency tables}
\author{Hidefumi Ohsugi \and Takayuki Hibi}
\date{}

\begin{document}

\maketitle

\begin{abstract}
In this paper, it is proved that,
if a toric ideal possesses a fundamental binomial
none of whose monomials is squarefree,
then the corresponding semigroup ring is not very ample.
Moreover, very ample semigroup rings of Lawrence type
are discussed.
As an application, we study very ampleness of configurations
arising from contingency tables.
\end{abstract}

\section{Introduction}
\label{intro}
A configuration in ${\mathbb R}^d$ is a finite set 
$A=\{ {\bf a}_1,\ldots,{\bf a}_n\} \subset {\mathbb Z}_{\geq 0}^d$
%which contained in a hyperplane in ${\mathbb R}^d$ without the origin.
such that there exists a vector ${\bf w} \in \RR^d$ satisfying
${\bf w} \cdot {\bf a}_i = 1$ for all $i$.
Let $K[{\bf t}]=K[t_1,\ldots,t_d]$
denote the polynomial ring in $d$ variables over a field $K$.
We associate a configuration $A$ with the semigroup ring
$K[A] =K[{\bf t}^{\a_1},\ldots, {\bf t}^{\a_n}]$,
where ${\bf t}^{\bf a} = t_1^{a_1} \cdots t_d^{a_d}$ if ${\bf a} = (a_1 ,\ldots,a_d)$.
Let $K[{\bf x}]=K[x_1,\ldots,x_n]$ denote the polynomial ring
in $n$ variables
over $K$.
The toric ideal $I_A$ of $A$ is the kernel of the surjective
homomorphism $\pi : K[{\bf x}] \longrightarrow K[A]$
defined by setting $\pi(x_i)= {\bf t}^{{\bf a}_i}$
for $1 \leq i \leq n$.

We are interested in the following conditions:
\begin{itemize}
\item[(i)]
$A$ is {\em unimodular}, i.e.,
the initial ideal of $I_A$ is generated by
squarefree monomials with respect to any monomial order;
\item[(ii)]
$A$ is {\em compressed}, i.e.,
the initial ideal of $I_A$ is generated by
squarefree monomials with respect to any 
reverse lexicographic order;
\item[(iii)]
there exists a monomial order $<$ such that
the initial ideal of $I_A$ with respect
to $<$ is generated by squarefree monomials;
\item[(iv)]
$K[A]$ is {\em normal}, i.e., $\NN A = \ZZ A \cap {\mathbb Q}_{\geq 0} A$;
\item[(v)]
$K[A]$ is {\em very ample}, i.e., $(\ZZ A \cap {\mathbb Q}_{\geq 0} A) \setminus \NN A$ 
is a finite (or empty) set.
\end{itemize}
Then
(i) $\Longrightarrow $
(ii) $\Longrightarrow $
(iii) $\Longrightarrow $
(iv) $\Longrightarrow $
(v) 
holds and each of the converse of them is false
in general.
%The inclusion
%${\mathbb N} A \subset \ZZ A \cap {\mathbb Q}_{\geq 0} A$
%holds in general.
If $K[A]$ is not normal, then an element of
$(\ZZ A \cap {\mathbb Q}_{\geq 0} A) \setminus \NN A$
is called {\em hole}.
%, which obstructs the
%sequential importance sampling???

%A semigroup ring $K[A]$ is called {\em normal}
%if we have $\NN A = \ZZ A \cap \QQ_{\geq 0} A$.
%If $K[A]$ is not normal, then an element of
%$(\ZZ A \cap \QQ_{\geq 0} A) \setminus \NN A$ is called a {\em hole}.
%A semigroup ring $K[A]$ is called {\em very ample}
%if $(\ZZ A \cap \QQ_{\geq 0} A) \setminus \NN A$ is a finite (or empty) set.

Let $P_A$ denote the convex hull of $A$.
For a subset $B \subset A$,
$K[B]$ is called {\it combinatorial pure subring} (\cite{cpure}, \cite{geom}) of $K[A]$
if there exists a face $F$ of $P_A$ such that $B = A \cap F$.
For example, if $K[B] = K[A] \cap K[ t_{i_1},\ldots,t_{i_s} ]$, then
$K[B]$ is a combinatorial pure subring of $K[A]$.
(This is the original definition of a combinatorial pure subring
in \cite{cpure}.)
A binomial $f \in I_A$ is called {\it fundamental}
if there exists a combinatorial pure subring $K[B]$ of $K[A]$
such that
$I_B$ is generated by $f$.
In Section \ref{fundamental}, it will be proved that,
if $I_A$ possesses a fundamental binomial
none of whose monomials is squarefree,
then $K[A]$ is not very ample.
The {\em Lawrence lifting} $\Lambda(\A)$ of
the configuration $\A$
is the configuration arising from the matrix
$$\Lambda(\A) = 
\left(
\begin{array}{cc}
\A & {\bf 0}\\
I_n & I_n
\end{array}
\right),
$$
where $I_n$ is the $n \times n$ identity matrix and
${\bf 0}$ is the $d \times n$ zero matrix.
A configuration $\A$ is called {\em Lawrence type}
if there exists a configuration $B$
such that $\Lambda(B) = \A$.
In Section \ref{fundamental}, it will be proved that a configuration of Lawrence type
is very ample if and only if it is unimodular.

In Section \ref{contingencytables}, 
by using the results in Section \ref{fundamental}, we study very ample
configurations arising from 
no $n$-way interaction models
for
$r_1 \times r_2 \times \cdots \times r_n$ contingency tables,
where $r_1 \geq r_2 \geq \cdots \geq r_n \geq 2$.
Let $A_{r_1 r_2 \cdots r_n}$ be the set of vectors
$
\e_{i_2 i_3 \cdots i_n}^{(1)} \oplus 
\e_{i_1 i_3 \cdots i_n}^{(2)} \oplus 
\cdots
\oplus
\e_{i_1 i_2 \cdots i_{n-1}}^{(n)},
$
where each $i_k$ belongs to $[r_k]=\{1,2,\ldots,r_k\}$
and $ \e_{j_1 j_2 \cdots j_{n-1}}^{(k)}   $
is a unit coordinate vector of 
$ {\ZZ}^{d_k} $ with $d_k = \frac{\prod_{\ell=1}^n r_\ell}{r_k}$.
The toric ideal $I_{\A_{r_1 r_2 \cdots r_n}}$
is the kernel of the %surjective 
homomorphism
$$
\pi: K[\{x_{i_1 i_2 \cdots i_n} \ | \ i_k \in [r_k]\}] \longrightarrow 
K[\{
%t_{\ \cdot \ i_2 i_3 \cdots i_n}^{(1)},
t_{i_1 \cdots i_{k-1}  i_{k+1} \cdots i_n}^{(k)}
%\ldots,
%t_{i_1 i_2 \cdots i_{n-1} \ \cdot \ }^{(n)}
\ | \ k \in [n],i_k \in [r_k]\}]
$$
defined by
$
\pi(x_{i_1 i_2 \cdots i_{n}}) = 
t_{i_2 i_3 \cdots i_n}^{(1)} 
t_{i_1 i_3 \cdots i_n}^{(2)} 
\cdots
t_{i_1 i_2 \cdots i_{n-1}}^{(n)}.
$
The following is known:

\bigskip

%\noindent
\begin{tabular}{|c|c|}
\hline
$r_1 \times r_2$ or $r_1 \times r_2 \times 2 \times \cdots \times 2$   &  unimodular\\
\hline
$r_1 \times 3 \times 3$ & compressed, not unimodular\\
\hline
$4 \times 4 \times 3$
& normal, not compressed\\
\hline
$5 \times 5 \times 3$ or $5 \times 4 \times 3$
& not compressed (normality is unknown)\\
\hline
otherwise, i.e.,
& \\
$n \geq 4$ and $r_3 \geq 3$ & not normal\\
$n =3$ and $r_3 \geq 4$ & \\
$n =3$, $r_3 = 3$, $r_1 \geq 6$ and $r_2 \geq 4$ & \\
\hline
\end{tabular}

\bigskip

\noindent
By virtue of the results in Section \ref{fundamental}, we will prove that
configurations in ``otherwise" part are not very ample.

\section{Fundamental binomials}
\label{fundamental}

The following lemma plays an important role in 
the present paper.

\begin{Lemma}
\label{tuika}
Let $K[B]$ be a combinatorial pure subring of $K[A]$.
If $K[A]$ is normal (resp. very ample), then
$K[B]$ is normal (resp. very ample).
\end{Lemma}

\begin{proof}
Let $K[B]$ be a combinatorial pure subring of $K[A]$.
It is enough to show that
$
(\ZZ B \cap \QQ_{\geq 0}B) \setminus \ZZ_{\geq 0} B
\subset
(\ZZ A \cap \QQ_{\geq 0}A) \setminus \ZZ_{\geq 0} A
$.

Let 
$\alpha \in (\ZZ B \cap \QQ_{\geq 0}B) \setminus \ZZ_{\geq 0} B.$
Since $B$ is a subset of $A$, we have $\alpha \in \ZZ A \cap \QQ_{\geq 0}A$.
Suppose that $\alpha \in \ZZ_{\geq 0} A $.
Then $\alpha = \sum_{{\bf a} \in A} z_{\bf a} {\bf a}$ with $0 \leq z_{\bf a} \in \ZZ$.
Since $\alpha \notin \ZZ_{\geq 0} B$,
$0 < z_{\bf a}$ for some ${\bf a} \in A \setminus B$.
Moreover, since $\alpha \in \QQ_{\geq 0}B $, 
$\alpha = \sum_{{\bf a} \in B} q_{\bf a} {\bf a}$ with $0 \leq q_{\bf a} \in \QQ$.
Thus
$\alpha =  \sum_{{\bf a} \in A} z_{\bf a} {\bf a} =
 \sum_{{\bf a} \in B} q_{\bf a} {\bf a}
.$
Since $K[B]$ is a combinatorial pure subring of $K[A]$,
there exists a face $F$ of $P_A$ such that $B=A\cap F$.
Then there exist ${\bf v} \in \RR^d$ and $c \in \RR$ satisfying
$$F = P_A \cap\{\ {\bf b} \in \RR^d \  |  \ {\bf v} \cdot {\bf b} = c \ \},$$
$$P_A \subset \{\ {\bf b} \in \RR^d \  |  \ {\bf v} \cdot{\bf b} \leq  c\ \}.$$
Then
${\bf v} \cdot {\bf a} = c$ for all ${\bf a} \in B$ and
${\bf v} \cdot {\bf a} < c$ for all ${\bf a} \in A \setminus B$.
Hence ${\bf v} \cdot \alpha =  c \  \sum_{{\bf a} \in B} q_{\bf a} <  c \ \sum_{{\bf a} \in A} z_{\bf a}$.
Thus we have $c \neq 0$ and $\sum_{{\bf a} \in B} q_{\bf a}  \neq \ \sum_{{\bf a} \in A} z_{\bf a}$.
On the other hand, since $A$ is a configuration,
there exists a vector ${\bf w} \in \RR^d$ satisfying
${\bf w} \cdot {\bf a} = 1$ for all ${\bf a} \in A$.
Hence ${\bf w} \cdot \alpha =   \sum_{{\bf a} \in B} q_{\bf a} = \sum_{{\bf a} \in A} z_{\bf a}$.
This is a contradiction.
Thus
$\alpha \in (\ZZ A \cap \QQ_{\geq 0}A) \setminus \ZZ_{\geq 0} A$
as desired.
\end{proof}

It is known \cite[Lemma 3.1]{cpure} that

\begin{Proposition}
\label{oldlemma}
If $g = u-v \in K[{\bf x}]$ is a binomial 
such that neither $u$ nor $v$ is squarefree
and if $I_A = (g)$, then $K[A]$ is not normal.
\end{Proposition}

We extend Proposition \ref{oldlemma} as follows:

\begin{Lemma}
\label{newlemma}
If $g = u-v \in K[{\bf x}]$ is a binomial 
such that neither $u$ nor $v$ is squarefree
and if $I_A = (g)$, then $K[A]$ is not very ample.
\end{Lemma}

\begin{proof}
Let $g = x_1^2 u' - x_2^2 v'$.
Since $g$ is irreducible,
$u' \ (\neq 1)$ is not divided by $x_2$ and 
$v' \ (\neq 1)$ is not divided by $x_1$.
Since $\pi(x_1^2 u') = \pi(x_2^2 v')$,
we have 
$\sqrt{\pi(u' v')} = \frac{\pi(x_1 u')}{\pi(x_2)}.$
Let $x_k$ be a variable with $k \neq 1,2$.
Then the monomial $\pi(x_k^m) \sqrt{\pi(u' v')}$
belongs to the quotient field of $K[A]$ and is integral over $K[A]$
for all positive integer $m$.
Suppose that there exists a monomial $w$ such that $\pi(w) = \pi(x_k^m) \sqrt{\pi(u' v')}$.
It then follows that the binomial $g' = x_1 u' x_k^m - x_2 w$ belongs to $I_A$.
Since $I_A= (g)$ and $x_1 u' x_k^m$ is divided by neither $x_1^2 u'$ nor $x_2^2 v'$,
we have $g'=0$.
Hence $x_2$ must divide $u'$, a contradiction.
Thus $\pi(x_k^m) \sqrt{\pi(u' v')}$ is a hole for all $m$ and
$K[A]$ is not very ample.
\end{proof}

\begin{Theorem}
\label{funda}
If $I_A$ possesses a fundamental binomial $g = u-v$
such that neither $u$ nor $v$ is squarefree,
then $K[A]$ is not very ample.
\end{Theorem}

\begin{proof}
Since $g$ is fundamental, 
there exists a combinatorial pure subring $K[B]$ of $K[A]$
such that
$I_B$ is generated by $g$.
Thanks to Lemma \ref{newlemma},
$K[B]$ is not very ample.
Since $K[B]$ is a combinatorial pure subring of $K[A]$,
$K[A]$ is not very ample by Lemma \ref{tuika}.
\end{proof}

Thanks to Theorem \ref{funda} together with
the results in
\cite{cpure},
we extend \cite[Theorem 3.4]{cpure} as follows:

\begin{Corollary}
\label{lawrence}
Let $K[A]$ be a semigroup ring and 
let $K[\Lambda(A)]$
its Lawrence lifting.
Then, the following conditions are equivalent:
\begin{itemize}
\item[(i)]
$K[A]$ is unimodular;
\item[(ii)]
$K[\Lambda(A)]$ is unimodular;
\item[(iii)]
%$K[\Lambda(A)]$ is normal;
%\item[(iv)]
$K[\Lambda(A)]$ is very ample.
\end{itemize}
\end{Corollary}

\begin{proof}
First, (ii) $ \Rightarrow $ (iii) is well-known.
On the other hand, (i) $\Leftrightarrow $ (ii) is proved in \cite[Theorem 3.4]{cpure}.

In order to show (iii) $ \Rightarrow $ (i), suppose that
$K[A]$ is not unimodular.
Then, by the same argument in Proof of \cite[Theorem 3.4]{cpure},
$I_{\Lambda(A)}$ has a fundamental binomial $\overline{g}$ none of whose monomials
is squarefree.
Thanks to Theorem \ref{funda}, $K[\Lambda(A)]$ is not very ample as desired.
\end{proof}

\begin{Remark}
{\em
A binomial $f$ belonging to $I_\A$ is called 
{\em  indispensable} 
if,
for an arbitrary
system ${\mathcal F}$ of binomial generators
of $I_\A$,
either $f$ or $-f$ appears in ${\mathcal F}$.
In particular, every fundamental binomial is indispensable.
However, Theorem \ref{funda} is not true if
we replace ``fundamental" with ``indispensable."
Let $K[A]= K[t_2, t_1 t_2, t_1^3 t_2, t_1^4 t_2] \subset K[t_1,t_2]$.
Then $K[A]$ is very ample and $I_A$ is generated by
the set of indispensable binomials
$\{  
x_1 x_4 - x_2 x_3, \ 
x_2^3 - x_1^2 x_3, \ 
x_3^3 - x_2 x_4^2, \ 
x_1 x_3^2 - x_2^2 x_4 
 \}$.
(The toric ideal $I_A$ has no fundamental binomials.)
}
\end{Remark}

\section{Configurations arising from contingency tables}
\label{contingencytables}

Configurations in ``otherwise" part of

\bigskip

%\noindent
\begin{tabular}{|c|c|}
\hline
$r_1 \times r_2$ or $r_1 \times r_2 \times 2 \times \cdots \times 2$   &  unimodular\\
\hline
$r_1 \times 3 \times 3$ & compressed, not unimodular\\
\hline
$4 \times 4 \times 3$
& normal, not compressed\\
\hline
$5 \times 5 \times 3$ or $5 \times 4 \times 3$
& not compressed (normality is unknown)\\
\hline
otherwise, i.e.,
& \\
$n \geq 4$ and $r_3 \geq 3$ & not normal\\
$n =3$ and $r_3 \geq 4$ & \\
$n =3$, $r_3 = 3$, $r_1 \geq 6$ and $r_2 \geq 4$ & \\
\hline
\end{tabular}

\bigskip

\noindent
are
studied in \cite{OHconti} by using the notion of combinatorial pure subring
%\cite{cpure}
and indispensable binomials.
% appearing in \cite{AoTa1}.
For $6 \times 4 \times 3$ case,
non-normality is shown in
\cite{Vlach}
and it was proved \cite{HTY} that it is not very ample.
On the other hand,
compressed configurations are classified in \cite{sullivant}.
For $4 \times 4 \times 3$ case,
it was announced in \cite[P.87]{HTY}
that
Ruriko Yoshida verified that it is normal
by using the software NORMALIZ (\cite{BrIc}).

The basic facts on $A_{r_1 \cdots r_n}$ are
\cite[Proposition 3.1 and Proposition 3.2]{OHconti}:

\begin{Proposition}
\label{pro1}
The configuration
$A_{r_1 \cdots r_n 2}$ is the Lawrence lifting of $A_{r_1 \cdots r_n}$.
\end{Proposition}

\begin{Proposition}
\label{pro2}
Suppose that $A_{r_1 \cdots r_n}$ and $A_{s_1 \cdots s_n}$ satisfy
$s_i \leq r_i$ for all $1 \leq i \leq n$.
Then $K[A_{s_1 \cdots s_n} ]$ is a combinatorial pure subring of 
$K[A_{r_1 \cdots r_n} ]$.
\end{Proposition}

\begin{Theorem}
Work with the same notation as above.
Then, each configuration in ``otherwise" part is not very ample.
\end{Theorem}

\begin{proof}
Let $A$ be a configuration in ``otherwise" part.
Thanks to Proposition \ref{pro2},
$K[A]$ has
at least one of
$K[A_{4   4   4}]$, $K[A_{6  4 3}]$ and
$K[A_{3 3 3 2 \cdots 2}]$
as a combinatorial pure subring.
It is easy to check that
$I_{A_{444}}$ has a fundamental binomial
$$
x_{1 1 1}^2 x_{1 3 3} x_{1 4 4} x_{2 2 3} x_{2 2 
          4} x_{2 3 2} x_{2 4 2} x_{3 1 3} x_{3 2 2} x_{3 
          4 1} x_{4 1 4} x_{4 2 2} x_{4 3 1} \\
$$
$$-    x_{1 1 3} x_{1 1 4} x_{1 3 1} x_{1 4 
          1} x_{2 2 2}^2 x_{2 3 3} x_{2 4 4} x_{3 1 
          1} x_{3 2 3} x_{3 4 2} x_{4 1 1} x_{4 2 4} x_{4 
          3 2},
$$
and
$I_{A_{643}}$ has a fundamental binomial
$$
x_{111} x_{221} x_{331} x_{641} x_{212} x_{522} x_{432} x_{642} x_{413} x_{323} x_{633}^2 x_{143} x_{543}
$$
$$
-
x_{211} x_{321} x_{631} x_{141} x_{412} x_{222} x_{632} x_{542} x_{113} x_{523} x_{333} x_{433} x_{643}^2.
$$
Since none of the monomials appearing above is squarefree, both
$K[A_{444}]$ and $K[A_{643}]$ are not very ample by Theorem \ref{funda}.
%(Note that, in \cite{HTY}, it is proved that $K[A_{643}]$ is not very ample.)
Moreover, since $A_{333}$ is not unimodular,
$K[A_{3 3  3   2  \cdots  2}]$
is not very ample by Corollary \ref{lawrence}
together with Proposition \ref{pro1}.
Thus, $K[A]$ is not very ample by Lemma \ref{tuika}.
\end{proof}

We close the present paper with an interesting problem.

\begin{Problem}
{\em
Find natural classes of configurations appearing in statistics
which is not normal but very ample.
}
\end{Problem}

\bigskip

\noindent
Hidefumi Ohsugi\\
Department of Mathematics, College of Science, Rikkyo University\\
Tokyo 171-8501, Japan \\
{\tt ohsugi@rkmath.rikkyo.ac.jp}

\bigskip

\noindent
Takayuki Hibi\\
Department of Pure and Applied Mathematics,
Graduate School of Information Science and Technology,
Osaka University\\
Toyonaka, Osaka 560-0043, Japan\\
{\tt hibi@math.sci.osaka-u.ac.jp}

\end{document}